\documentclass[conf]{new-aiaa}
\usepackage[utf8]{inputenc}

\usepackage{graphicx}
\usepackage{amsmath}

\usepackage{amsthm}
\usepackage[version=4]{mhchem}
\usepackage{siunitx}
\usepackage{longtable,tabularx}
\usepackage{jhoagg}
\usepackage{cleveref}
\usepackage{subcaption}
\usepackage{tikz}
\usetikzlibrary{arrows.meta, shapes.geometric, positioning}

\newtheorem{theorem}{Theorem}

\newtheorem{lemma}{Lemma}

\newtheorem{definition}{Definition}
\newtheorem{fact}{Fact}

\newtheorem{problem}{Problem}

\newtheorem{example}{Example}

\allowdisplaybreaks

\newcounter{enumA}
\newcounter{enumC}
\newcounter{enumO}

\newlist{enumA}{enumerate}{1}
\setlist[enumA,1]{%
  label=(A\arabic*),
  leftmargin=1cm,
  before=\setcounter{enumAi}{\value{enumA}},
  after=\setcounter{enumA}{\value{enumAi}}
}

\newlist{enumC}{enumerate}{1}
\setlist[enumC,1]{%
  label=(C\arabic*),
  leftmargin=1cm,
  before=\setcounter{enumCi}{\value{enumC}},
  after=\setcounter{enumC}{\value{enumCi}}
}

\newlist{enumO}{enumerate}{1}
\setlist[enumO,1]{%
  label=(O\arabic*),
  leftmargin=1cm,
  before=\setcounter{enumOi}{\value{enumO}},
  after=\setcounter{enumO}{\value{enumOi}}
}

\title{Geometric Conditions for Lossless Convexification in Linear Optimal Control with Discrete-Valued Inputs: Real-Time Implementation for Spacecraft Rendezvous}

\author{Felipe Arenas-Uribe\footnote{Graduate Research Assistant, Mechanical and Aerospace Engineering, AIAA Student Member.}, Hasan A. Poonawala\footnote{Associate Professor, Mechanical and Aerospace Engineering.} and Jesse B. Hoagg\footnote{Professor, Mechanical and Aerospace Engineering, Associate Fellow AIAA.}}
\affil{University of Kentucky, Lexington, Kentucky, 40506}

\begin{document}

\maketitle

\begin{abstract}
    Optimal control problems with discrete-valued inputs are inherently challenging due to their mixed-integer nature, rendering them generally intractable for real-time, safety-critical aerospace applications. 
    Lossless convexification offers a powerful alternative by reformulating these mixed-integer programs into computationally efficient convex programs. 
    This paper develops a lossless convexification framework for the optimal control of linear time-varying systems with discrete-valued inputs. 
    We extend existing theoretical results by demonstrating that system normality is preserved when reformulating Lagrange-form problems into Mayer-form via an epigraph transformation. 
    Furthermore, we establish that under simple geometric conditions on the input set, the solution to the relaxed convex problem strictly satisfies the original non-convex input constraints. 
    This framework enables the real-time computation of optimal discrete-valued controls without resorting to mixed-integer optimization. 
    The proposed algorithm is validated on a spacecraft rendezvous maneuver utilizing discrete-valued reaction thrusters in an elliptical orbit. 
    Numerical results from Monte Carlo simulations show that the finite-grid algorithm consistently yields near-discrete control inputs, with computational times compatible with safety-critical, on-board applications.
\end{abstract}

\section{Introduction}

Discrete-valued inputs are common in physical systems, as many actuators operate in an on–off manner or with finitely many magnitude levels. 
% Spacecraft performing proximity maneuvers such as rendezvous and capture have actuator configurations that predominantly rely on discrete-valued actuators such as reaction control thrusters \cite{pasand_study_2017}. 
Spacecraft performing proximity maneuvers such as rendezvous and capture have actuator configurations that predominantly rely on discrete-valued actuators such as reaction control thrusters. 
Optimal control of such systems is computationally challenging since the discrete actuator behavior introduces mixed-integer constraints. 
Mixed-Integer Programs (MIPs) can numerically solve such control problems, but are inherently difficult to solve, lacking general convergence guarantees and exhibiting exponential growth in runtime with the number of decision variables \cite{junger_mixed_2013}. 
These computational challenges render MIPs unsuitable for real-time applications such as autonomous vehicle control in safety-critical applications.

Recent work has used learning-based methods to warm-start MIP solvers and thereby accelerate their convergence \cite{cauligi_coco_2022, cauligi_learning_2020}. These methods use neural networks to provide an initial, potentially suboptimal solution that an optimizer can refine. 
Nonetheless, these strategies provide no formal guarantees of optimality or convergence and require a pre-computed dataset of solutions to train the neural network. 
An alternative to MIP solvers is the use of homotopy methods. 
Their core idea is to introduce activation functions and sequentially solve a continuous problem, adjusting the activation parameters such that the solution takes a mixed-integer behavior \cite{saranathan_relaxed_2018,malyuta_fast_2023,arya_composite_2021}. 
However, homotopy methods introduce approximations, can be computationally expensive, and do not guarantee optimality. 
In contrast, lossless convexification methods have shown that certain MIPs can be reformulated as equivalent convex programs (CPs) without loss of optimality, offering a computationally tractable alternative \cite{malyuta_convex_2022}. 
A lossless convexification for optimal control problems (OCP) in Mayer-form with linear time-invariant (LTI) systems and disconnected input sets was introduced in \cite{harris_optimal_2021}, where system normality and an extreme-point relaxation were exploited to achieve an exact convex reformulation. 
Similarly, \cite{woodford_geometric_2022} exploited strong observability to obtain a lossless convexification for a similar problem with a family of state constraints. 
It was proposed in \cite{harris_optimal_2021} that Lagrange-form problems can be reformulated as Mayer-form problems via an epigraph transformation; however, it remains unclear whether system normality is preserved under this transformation. 

Lagrange-form OCPs include the fuel-optimal control problem, where the $L^1$-norm of the control input serves as the cost function since it is proportional to the total control effort \cite{ross_how_2004}. 
This problem is particularly relevant in spacecraft applications, where mission lifetime is limited by onboard propellant and maneuvers must therefore be fuel-optimal. 
Standard trajectory generation algorithms for spacecraft rely in quadratic cost functions which don't minimize fuel consumption directly, thus reducing mission lifetime \cite{ross_space_2006}. 
In addition, spacecraft employ discrete-valued actuators such as cold-gas thrusters in reaction control systems \cite{pasand_study_2017,weston_state---art_2025}. Hence, it is essential to develop real-time fuel-optimal control strategies that explicitly account for discrete inputs.

\begin{figure}[t]
\centering
\includegraphics[width=0.55\linewidth]{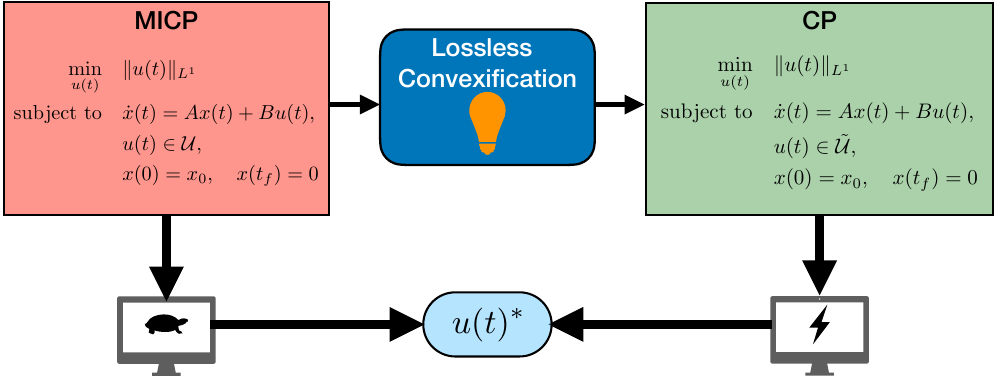}
%\vspace{-0.15 in}
\caption{Lossless convexification enables real-time solution of optimal control problems with discrete-valued inputs by solving a convex program.}
\label{fig:Framework}
\vspace{-0.15 in}
\end{figure}

In this paper, we extend lossless convexification results from \cite{harris_optimal_2021} and show that Lagrange-form OCPs for linear systems with discrete-valued inputs can be reformulated as a Mayer-form MIP while preserving system normality. 
We demonstrate that the resulting Mayer-form MIP can be solved as a convex program by leveraging system normality and convex hull relaxation of the input set (see \Cref{fig:Framework}), granted that the input set satisfies a geometric condition.
To bridge theory and practice, we present a finite-dimensional formulation of the convex OCP alongside a bisection algorithm designed to identify a feasible time horizon. 
The proposed framework is validated through numerical experiments on a spacecraft rendezvous maneuver in an elliptical orbit utilizing discrete-valued reaction thrusters. 
Finally, we evaluate the impact of discretization coarseness on solution discreteness and conduct Monte Carlo simulations across diverse initial conditions, showing that the finite-grid algorithm consistently yields near-discrete controls with computation times compatible with real-time, safety-critical applications.

The remainder of this paper is organized as follows. \Cref{sec:ProblemForm} outlines the problem formulation. \Cref{sec:EquivFormulation} presents an equivalent formulation of the problem. \Cref{sec:LosslessConv} presents the corresponding lossless convexification. Numerical simulation results are provided in \Cref{sec:NumericalSim}, and \Cref{sec:Conclusions} concludes with final remarks and directions for future research.

\section{Problem Formulation}\label{sec:ProblemForm}

\subsection{Notation}

Let $\mathbb{N}_{>0}$ be the set of positive natural numbers not including 0. Let $\BBR_{\geq 0}$ be the set of non-negative real numbers. Let $\SA \subset \mathbb{R}^n$. Then, the convex hull of $\SA$ is $\mathrm{conv}(\SA)$, the extreme points of $\SA$ are $\mathrm{ex}(\SA)$, and the boundary of $\SA$ is $\mathrm{bd}(\SA)$. The vector $e_i$ denotes the $i$-th basis vector of $\mathbb R^m$. Let $f:\mathbb{R}^{n} \to \mathbb{R}$, the \emph{epigraph} of $f$ is $\mathrm{epi}(f) \triangleq \{(x,\mu) \in \mathbb{R}^{n+1}: f(x) \leq \mu\}$. 

\subsection{Discrete-Valued Optimal Control Problem}

Consider the linear time-varying system
\begin{equation}\label{eq:linear_system}
    \dot{x}(t) = A(t) x(t) + B(t) u(t),
\end{equation}
where $A(t) \in \mathbb{R}^{n \times n}$, $B(t) \in \mathbb{R}^{n \times m}$, the pair $\big(A(t),B(t)\big)$ is controllable for all $t\geq0$, $x(t) \in \mathbb{R}^n$ is the state, $x(0) = x_0 \in \mathbb{R}^n \setminus \{0\}$ is the initial condition, and $u(t) \in \mathcal{U} \subset \mathbb{R}^m$ is the input. The input set $\mathcal{U}$ is a finite disconnected set
\begin{equation}\label{eq:input_set}
    \mathcal{U} \triangleq \{u_1, u_2, \ldots, u_{\ell}\},
\end{equation}
where $u_i \in \mathbb{R}^m,\; i=1,\ldots,\ell$ and $\ell \in \mathbb{N}_{>0}$. 

Consider the cost function $\psi: \mathbb{R}^m \to \BBR_{\geq 0}$, where $\psi(u)$ is a positive-definite polyhedral convex function. 
Note that $\psi(u)$ is polyhedral convex if and only if $\mathrm{epi}(\psi)$ is a polyhedral convex set \cite[Corollary 19.2.1]{rockafellar_convex_2015}.
Let $\mu \in [0, \bar{u}]$, where $\bar{u} \triangleq \max_{u \in \SU} \psi(u)$. Then, define the set
\begin{equation}\label{eq:relaxed_set}
    \tilde{\SU}_e \triangleq \mathrm{conv}(\SU) \times [0, \bar{u}] \cap \mathrm{epi}(\psi).
\end{equation}
We make the assumption:
\begin{enumA}
    \item $\mathrm{ex}(\tilde{\SU}_e) = \{(u,\mu) \in \mathbb{R}^{m+1}: u \in \SU, \mu = \psi(u)\}.$ \label{ass:augmented_set_vertices}
\end{enumA}
This assumption implies that the extreme points of the set $\tilde{\SU}_e$ are members of the input set $\SU$. This geometric property is characterized by $\psi$ and $\SU$. This property will be useful in later sections to enable the lossless convexification of an OCP. The following examples show pairs of cost functions and sets which satisfy \ref{ass:augmented_set_vertices}.

\begin{figure}[h]
    \centering
    % First Subfigure (Left)
    \begin{subfigure}[b]{0.48\linewidth}
        \centering
        \includegraphics[width=0.9\linewidth]{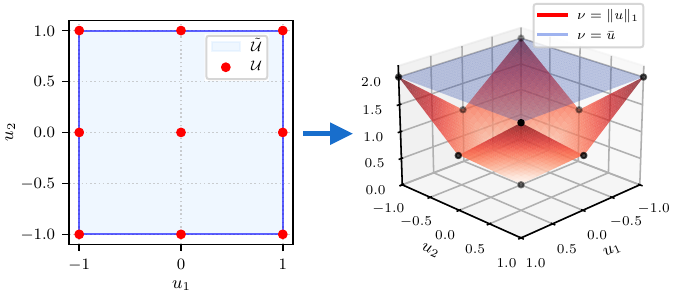}
        \caption{Input set $\SU$ and its augmented set $\tilde{\mathcal{U}}_e$ from \Cref{ex:set_example_1}.}
        \label{fig:example1_polytope}
    \end{subfigure}
    \hfill % Adds horizontal space to push the subfigures to the edges
    % Second Subfigure (Right)
    \begin{subfigure}[b]{0.48\linewidth}
        \centering
        \includegraphics[width=0.9\linewidth]{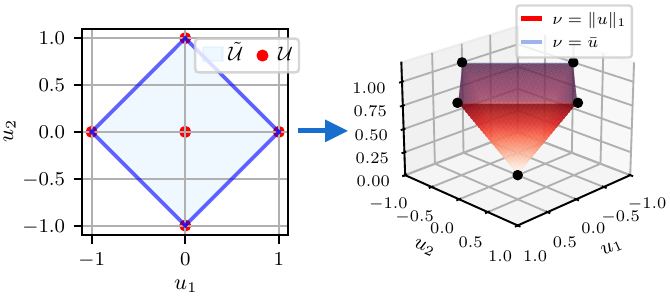}
        \caption{Input set $\SU$ and its augmented set $\tilde{\mathcal{U}}_e$ from \Cref{ex:set_example_2}.}
        \label{fig:example2_polytope}
    \end{subfigure}
    
    \caption{Visualization of the sets $\SU$ and $\tilde \SU_e$ from \Cref{ex:set_example_1,ex:set_example_2}.}
    \label{fig:combined_polytopes} % Main figure label for referencing both together
\end{figure}

\begin{example}\label{ex:set_example_1}\rm
    Let $\psi(u) = \|u\|_1$ and consider a system whose input channels can be actuated independently but are limited between two discrete values. Let $u_{\max} > 0$ and consider the input set
    \begin{equation*}
        \mathcal{U} \triangleq \{-u_{\max},0,u_{\max}\}^m.
    \end{equation*}
    As seen in \Cref{fig:example1_polytope}, \ref{ass:augmented_set_vertices} is satisfied since all vertices of $\tilde{\SU}_e$ are elements of the set $\SU$. This type of actuator set is common in systems with discrete-valued actuators (e.g. valves or reaction thrusters) who operate independently of each other.
\end{example}

\begin{example}\label{ex:set_example_2}\rm
    Let $\psi(u) = \|u\|_1$ and consider a system whose input channels are limited by a magnitude $u_{\max}>0$, and the admissible combinations are constrained such that the 1-norm of the input vector does not exceed this bound. Let $u_{\max} > 0$ and consider the input set
    \begin{equation*}
        \mathcal{U} \triangleq \{0\}\cup\{\pm u_{\max} e_i : i = 1,\ldots,m\}.
    \end{equation*}
    This structure, shown in \Cref{fig:example2_polytope}, satisfies \ref{ass:augmented_set_vertices} since all vertices of $\tilde{\SU}_e$ are elements of the set $\SU$. The input set considered in this example is common in systems where multiple actuators share a limited total control authority, such that at most one actuator can be fully engaged at any given time. This is the case in small satellite reaction control systems which share a singular propellant tank \cite{weston_state---art_2025}.
\end{example}

Let $t_f>0$ be the time horizon. The objective is to solve the optimal control problem:
\begin{problem}\rm{Lagrange-Form Mixed-Integer Problem}\label{prob:MICP}
    \begin{subequations}
        \begin{align}
            \min \quad & \int_{0}^{t_f} \psi(u(t)) \, dt \label{eq:cost} \\
            \text{subject to} \quad 
            & \dot{x}(t) = A(t) x(t) + B(t) u(t),\\
            & u(t) \in \mathcal{U}, \label{eq:input_constraints} \\
            & x(0) = x_0, \quad x(t_f) = 0. \label{eq:boundary_conditions}
        \end{align}
    \end{subequations}
\end{problem}
We assume that the time horizon is selected such that there exists $u:[0,t_f] \to \SU$ for which the solution of \cref{eq:linear_system} satisfies $x(t_f)=0$. 
This assumption ensures that the feasible set of \Cref{prob:MICP} is nonempty.

\Cref{prob:MICP} offers a general framework capable of capturing a diverse class of OCPs with relevance to aerospace systems. 
For instance, setting $\psi(u) = \|u\|_1$ yields a fuel-optimal control formulation, which is central to spacecraft maneuvers such as satellite rendezvous, proximity operations, docking, and station-keeping. 
On the other hand, setting \(\psi(u)=\|u\|_\infty\) yields an accumulated peak-channel control formulation, which penalizes the largest input magnitude at each instant.
Alternatively, choosing the zero-pseudonorm $\psi(u) = \|u\|_0$ yields the maximum hands-off control problem \cite{nagahara_maximum_2016}. 
This formulation is pertinent to scenarios requiring sparse actuator intervention, finding applications in orbital station-keeping and trajectory maintenance.

\section{Equivalent Formulation in Mayer-form}\label{sec:EquivFormulation}

To analyze the convexification properties of \Cref{prob:MICP}, we reformulate the original Lagrange-form OCP into an equivalent Mayer-form OCP. 
By invoking an epigraph transformation, we embed the integral cost into the state \cite{boyd_convex_2023}. This allows us to leverage existing convexification frameworks in optimal control literature, such as those established in \cite{harris_optimal_2021}. First, define the augmented input set
\begin{equation}\label{eq:def_augmented_set}
    \mathcal{U}_e \triangleq \{(u,\mu) \in \mathbb{R}^{m+1} : u \in \mathcal{U}, \mu \in [0, \bar{u}],\; \psi(u) \leq \mu\}.
\end{equation}
Next, let $x_c(t) \in \BBR$ be the cost state which satisfies
\begin{equation}\label{eq:cost_state_dynamics}
    \dot{x}_c(t) = \mu(t),
\end{equation}
with the initial condition $x_c(0) = 0$. Define $\hat{u} \triangleq \begin{bmatrix} u & \mu \end{bmatrix}$. Then, the cascade of \cref{eq:linear_system} and \cref{eq:cost_state_dynamics} is
\begin{equation} \label{eq:augmented_dynamics}
    \dot{\hat{x}}(t) = \hat{A}(t) \hat{x}(t) + \hat{B}(t)\hat{u}(t),
\end{equation}
where
\begin{equation*}
    \hat{x} \triangleq \begin{bmatrix} x \\ x_c \end{bmatrix}, \quad
    \hat{A}(t) \triangleq \begin{bmatrix}
        A(t) & 0 \\
        0 & 0
    \end{bmatrix}, \quad \hat{B}(t) \triangleq \begin{bmatrix}
        B(t) & 0 \\
        0 & 1
    \end{bmatrix}.
\end{equation*}
Thus, consider the Mayer-form OCP:
\begin{problem}\rm{Mayer-form Mixed-Integer Problem}\label{prob:MICP_Mayer}
    \begin{subequations}
        \begin{align}
            \min \quad & x_c(t_f) \label{eq:TerminalCost}\\
            \text{subject to} \quad 
            & \dot{\hat{x}}(t) = \hat{A}(t) \hat{x}(t) + \hat{B}(t)\hat{u}(t),\\
            & \hat u \in {\mathcal{U}}_e, \label{eq:input_cons}\\
            & x(0) = x_0, \quad x(t_f) = 0, \\
            & x_c(0) = 0.
        \end{align}
    \end{subequations}
\end{problem}

We now study the equivalency of \Cref{prob:MICP_Mayer} with respect to \Cref{prob:MICP} from an optimal control perspective. Define the Hamiltonian of \Cref{prob:MICP}
\begin{equation}\label{eq:Hamiltonian_1}
    H_1(t,x,u,p) \triangleq \psi(u(t)) + p(t)^\top \big( A(t)x(t)+B(t)u(t)\big),
\end{equation}
where $p(t) \in \BBR^n$ is the co-state associated with \cref{eq:linear_system}. The co-state satisfies
\begin{equation}
    \dot p(t) = - \frac{H_1(t,x,u,p)}{\partial x} = -A(t)^\top p(t),
\end{equation}
with the initial condition $p(0) = p_0 \in \BBR^n \setminus \{0\}$ given by the transversality condition \cite[Theorem 6.3.5]{berkovitz_nonlinear_2013}. It follows from the necessary conditions of optimal control \cite[Theorem 6.7.1]{berkovitz_nonlinear_2013}, that the optimal control of \Cref{prob:MICP} $u^*:[0,t_f] \to \SU$ satisfies
\begin{equation}\label{eq:optimal_control_p1}
    u^*(t) = \arg \min_{u \in \SU} \big[ H_1(t,x,u,p) \big].
\end{equation}

Next, define the Hamiltonian of \Cref{prob:MICP_Mayer}
\begin{equation}\label{eq:Hamiltonian_2}
    H_2(t,\hat x,\hat u,p, p_c) \triangleq p(t)^\top \big( A(t)x(t)+B(t)u(t)\big) + p_c(t) \mu(t),
\end{equation}
where $p_c(t) \in \BBR$ is the co-state associated with \cref{eq:cost_state_dynamics}. 
The terminal cost of \Cref{prob:MICP_Mayer} yields the transversality condition $p_c(t_f)=\frac{\partial x_c(t_f)}{\partial x_c}=1$. 
Since $\dot p_c(t) = -\partial H_2/\partial x_c = 0$, the co-state $p_c$ is constant over $[0,t_f]$, thus for all $t\in[0,t_f]$, $p_c(t) = 1$. 
Next, it follows from the necessary conditions of optimal control \cite[Theorem 6.7.1]{berkovitz_nonlinear_2013}, that the optimal control of \Cref{prob:MICP_Mayer} $\hat u^*: [0,t_f] \to \SU_e$ satisfies 
\begin{equation}\label{eq:optimal_control_p2}
    \hat u^*(t) = \arg \min_{\hat u \in \SU_e} \big[ H_2(t,\hat x,\hat u,p, p_c) \big].
\end{equation}

The following result shows that the optimal control of \Cref{prob:MICP_Mayer} is also the optimal control of \Cref{prob:MICP}.

\begin{lemma}\rm\label{lemma:MICP2mayer}
Let $u_\star:[0,t_f] \to \SU$ and $\mu_\star:[0,t_f] \to \BBR$ such that for all $t \in [0,t_f]$, $\big(u_\star(t),\mu_\star(t)\big) = \hat{u}^*(t)$. Then, for all $t \in [0,t_f]$, $\mu_\star(t) = \psi(u_\star(t))$ and $u_\star(t)= u^*(t)$.
\end{lemma}

\begin{proof}
Let $t\in[0,t_f]$. 
It follows from \cref{eq:optimal_control_p2} that $\bigl(u_\star(t),\mu_\star(t)\bigr)\in\SU_e$, and $H_2\bigl(t,\hat x,(u_\star(t),\mu_\star(t)),p,p_c\bigr) \leq H_2\bigl(t,\hat x, (u_\star(t),\psi(u_\star(t))),p,p_c\bigr)$.
Furthermore, it follows from \cref{eq:def_augmented_set} that $\mu_\star(t)\geq\psi(u_\star(t))$.
Since $p_c(t)=1$, \cref{eq:Hamiltonian_1,eq:Hamiltonian_2} yields
\begin{equation*}
    H_1(t,x,u_\star(t),p) + \mu_\star(t)-\psi(u_\star(t)) \leq H_1(t,x,u_\star(t),p).
\end{equation*}
Therefore, $\mu_\star(t)\leq\psi(u_\star(t))$. Moreover, it follows from \cref{eq:relaxed_set} that $\mu_\star(t)\geq\psi(u_\star(t))$. Thus, $\mu_\star(t)=\psi(u_\star(t))$.

Next, let $u_1(t)\in\SU$ such that $u_1(t) \neq u_\star(t)$ and let $\mu_1(t) = \psi(u_1(t))$. 
Note that $(u_1(t), \mu_1(t)) \in \SU_e$.
It follows from \cref{eq:optimal_control_p2} that $H_2\bigl(t,\hat x,(u_\star(t),\psi(u_\star(t))),p,p_c\bigr) \leq H_2\bigl(t,\hat x,(u(t),\psi(u(t))),p,p_c\bigr)$.
Since, in addition, $p_c(t)=1$, \cref{eq:Hamiltonian_1,eq:Hamiltonian_2} yields
\begin{equation*}
    H_1(t,x,u_\star(t),p) \leq H_1(t,x,u(t),p).
\end{equation*}
Thus, $u_\star(t) = \arg\min_{u\in\SU}H_1(t,x,u,p)$. Hence, it follows from \cref{eq:optimal_control_p1} that $u_\star(t) = u^*(t)$.
\end{proof}

\Cref{lemma:MICP2mayer} implies that the optimal control of \Cref{prob:MICP_Mayer} is given by the pair $\hat u^*(t) = \big ( u^*(t), \mu^*(t)\big)$, where the first element is the optimal control of \Cref{prob:MICP}. Thus, the next section focuses on developing a tractable solution to \Cref{prob:MICP_Mayer}.

\section{Normality for Lossless Convexification}\label{sec:LosslessConv}

In this section, we develop a lossless convexification of \Cref{prob:MICP_Mayer}.
Consider the relaxed augmented set $\tilde{\mathcal{U}}_e$ in \cref{eq:relaxed_set}, which relaxes the input constraint \cref{eq:input_cons} by using a convex hull relaxation of the original input set (i.e. $u \in \mathrm{conv}(\mathcal{U})$). 
The following result shows that by using the set $\tilde{\mathcal{U}}_e$ as the input constraint, we introduce a convex relaxation of \cref{eq:input_cons}.

\begin{fact}\label{fact:augmented_set_polytope}
The set $\tilde{\mathcal{U}}_e$ is a polytope.
\end{fact}

\begin{proof}
Note that $\mathrm{conv}(\SU)$ is the convex hull of a finite set of points, thus it is a polytope \cite[Theorem 18.5]{rockafellar_convex_2015}.
Next, note that $[0, \bar{u}]$ is a one-dimensional simplex, which is a polytope by definition.
Furthermore, since $\psi(u)$ is a polyhedral convex function, it follows from \cite[Corollary 19.2.1]{rockafellar_convex_2015} that $\mathrm{epi}(\psi)$ is a polyhedral convex set. 
Since, in addition, $\mathrm{conv}(\SU)$ and $[0, \bar{u}]$ are bounded, then it follows from \cref{eq:relaxed_set} that $\tilde{\mathcal{U}}_e$ is a polytope.
\end{proof}

% Since the set $\tilde{\mathcal{U}}_e$ is a polytope, then we can use its half-space representation to enforce the input constraints as linear constraints. 
Then, consider the convex relaxation of \Cref{prob:MICP_Mayer}:
\begin{problem}\rm{Mayer-form Convex Problem}\label{prob:CP_Mayer}
    \begin{subequations}
        \begin{align}
            \min \quad & x_c(t_f)\\
            \text{subject to} \quad 
            & \dot{\hat{x}}(t) = \hat{A}(t) \hat{x}(t) + \hat{B}(t) \hat{u}(t), \label{eq:Prob_Aug_Dynamics}\\
            & \hat{u} \in \tilde{\mathcal{U}}_e,\\
            & x(0) = x_0, \quad x(t_f) = 0, \\
            & x_c(0) = 0.
        \end{align}
    \end{subequations}
\end{problem}
It follows from the necessary conditions of optimal control \cite[Theorem 6.7.1]{berkovitz_nonlinear_2013} that the optimal control of \Cref{prob:CP_Mayer} $\tilde u^*: [0,t_f] \to \SU_e$ satisfies
\begin{equation}\label{eq:optimal_control_p3}
    \tilde u^*(t) = \arg \min_{\tilde u \in \tilde{\SU}_e} \big[ H_2(t,\hat x,\tilde u,p, p_c) \big].
\end{equation}

If the minimizer of \cref{eq:optimal_control_p3} over $\tilde{\mathcal{U}}_e$ is not unique, then \cref{eq:optimal_control_p3} does not determine a well-defined control law. To exclude this degeneracy, the concept of normality is defined \cite[Definition 6.7.11]{berkovitz_nonlinear_2013}.

\begin{definition}\rm \label{def:normality}
    The system $(\hat{A},\hat{B})$ is \emph{strongly normal} with respect to $\tilde{\mathcal{U}}_e$ if, for every $p_0 \in \BBR^{n} \setminus \{0\}$ and for all but a finite set of points $t \in [0,t_f]$, \eqref{eq:optimal_control_p3} has a unique solution $\hat{u}^*(t) \in \tilde{\mathcal{U}}_e$.
\end{definition}

A fundamental result of normality for linear systems with polytopic input sets is presented next \cite[Corollary 6.7.7]{berkovitz_nonlinear_2013}.

\begin{lemma}\label{lemma:vertices_normality}
    Assume that $(\hat{A},\hat{B})$ is strongly normal with respect to $\tilde{\mathcal{U}}_e$.
    Then, $\tilde u^*(t)$ is piecewise constant on $[0,t_f]$ and for almost every $t \in [0,t_f]$, $\tilde{u}^*(t) \in \mathrm{ex}(\tilde{\mathcal{U}}_e)$.
\end{lemma}

This result implies that the optimal control $\tilde{u}^*$ can only take values from a finite set given by the vertices of $\tilde{\mathcal{U}}_e$. The next result studies normality of the original system and the augmented system through the epigraph transformation.

\begin{lemma}\rm{ \label{lemma:Normality_augmented}}
    % $(A, B)$ is strongly normal with respect to $\mathrm{conv}(\SU)$ if and only if $(\hat{A}, \hat{B})$ is strongly normal with respect to $\tilde{\mathcal{U}}_e$.
    $(A, B)$ is strongly normal with respect to $\mathrm{conv}(\SU)$ if and only if $(\hat{A}, \hat{B})$ is strongly normal with respect to $\tilde{\mathcal{U}}_e$.
\end{lemma}

\begin{proof}
    To prove sufficiency, suppose that $(A,B)$ is strongly normal with respect to $\mathrm{conv}(\SU)$. Then, for every $p_0\in\BBR^n\setminus\{0\}$ and for almost every $t \in [0,t_f]$,
    \begin{equation}\label{eq:u_1_opt}
        u_1^*(t) = \arg\min_{u\in\mathrm{conv}(\SU)} \big[ p(t)^\top\big(A(t)x(t)+B(t)u(t)\big) +\psi(u) \big],
    \end{equation}
    has a unique solution.

    Let $\tilde u^*(t)=(u^*(t),\mu^*(t))$ be a minimizer of \cref{eq:optimal_control_p3}. Following the same argument used in \Cref{lemma:MICP2mayer}, with $\SU$ replaced by $\mathrm{conv}(\SU)$, it follows that $\mu^*(t)=\psi(u^*(t))$ and that $u^*(t) = \arg\min_{u\in\mathrm{conv}(\SU)} \big[ p(t)^\top\big(A(t)x(t)+B(t)u(t)\big) +\psi(u) \big]$.
    It follows from \cref{eq:u_1_opt} that $u^*(t)=u_1^*(t)$. Therefore, the minimizer of \cref{eq:optimal_control_p3} satisfies $\tilde u^*(t) = \big(u_1^*(t),\psi(u_1^*(t))\big)$.
    Hence, $\tilde u^*(t)$ is unique, which implies that $(\hat A,\hat B)$ is strongly normal with respect to $\tilde{\mathcal{U}}_e$.

    Conversely, suppose that $(\hat A,\hat B)$ is strongly normal with respect to $\tilde{\mathcal{U}}_e$. Then, for every $p_0\in\BBR^n\setminus\{0\}$ and for almost every $t\in[0,t_f]$, \cref{eq:optimal_control_p3} has a unique solution $\tilde u^*(t) = \big(u^*(t),\psi(u^*(t))\big)$.
    Next, since $p_c(t)=1$, it follows from \cref{eq:Hamiltonian_1,eq:Hamiltonian_2} that $u^*(t)$ satisfies $u^*(t) = \arg\min_{u\in\mathrm{conv}(\SU)} \big[ p(t)^\top\big(A(t)x(t)+B(t)u(t)\big) +\psi(u(t)) \big]$.
    Suppose, for contradiction, that there exists $u_2(t) \in \mathrm{conv}(\SU)$ such that $u_2^*(t) = \arg\min_{u\in\mathrm{conv}(\SU)} \big[ p(t)^\top\big(A(t)x(t)+B(t)u(t)\big) +\psi(u(t)) \big]$ and $u_2^*(t)\neq u^*(t)$. Then, $\big(u_2^*(t),\psi(u_2^*(t))\big) \in\tilde{\mathcal{U}}_e$, would also be a minimizer of \cref{eq:optimal_control_p3}. This contradicts the uniqueness of $\tilde u^*(t)$. Therefore, $u^*(t)$ is unique, which implies that $(A,B)$ is strongly normal with respect to $\mathrm{conv}(\SU)$.
\end{proof}

This result is useful since it allows to analyze the normality of the augmented system \cref{eq:augmented_dynamics} with respect to $\tilde{\SU}_e$ based on the lower-dimensional system \cref{eq:linear_system} and the set $\mathrm{conv}(\SU)$. Note that the relaxed input set $\mathrm{conv}(\SU)$ can be easily computed offline using a convex hull algorithm, such as the Double-Description method \cite{fukuda_double_1996}. This algorithm exhibits super-polynomial complexity in the worst case, which does not impede real-time performance, as the half-space representation of $\mathrm{conv}(\SU)$ only needs to be computed once offline. The algorithm only requires re-execution if the actuator input set $\mathcal{U}$ is modified.

The next result is the main result, which shows that the optimal control of \Cref{prob:CP_Mayer} is also the optimal control of \Cref{prob:MICP}.

\begin{theorem}\label{thm:CP_for_MICP}
Assume that $(A,B)$ is strongly normal with respect to $\mathrm{conv}(\SU)$ and that \ref{ass:augmented_set_vertices} is satisfied. 
Let $u_\star:[0,t_f] \to \mathrm{conv}(\SU)$ and $\mu_\star:[0,t_f] \to \BBR$ such that for all $t \in [0,t_f]$, $\big(u_\star(t),\mu_\star(t)\big) = \tilde{u}^*(t)$, where $\tilde{u}^*(t)$ is given by \cref{eq:optimal_control_p3}.
Then, $u_\star$ is an optimal control of \Cref{prob:MICP}. Furthermore, $\big(u_\star(t),\mu_\star(t)\big)$ is piecewise constant on $[0,t_f]$ and, for almost every $t\in[0,t_f]$, $u_\star(t)\in\SU$.
\end{theorem}

\begin{proof}
By \Cref{lemma:Normality_augmented}, the strong normality of $(A,B)$ with respect to $\mathrm{conv}(\SU)$ implies that $(\hat{A},\hat{B})$ is strongly normal with respect to $\tilde{\mathcal{U}}_e$. Furthermore, \ref{ass:augmented_set_vertices} implies that $\mathrm{ex}(\tilde{\mathcal{U}}_e) \subseteq \mathcal{U}_e$.
Therefore, it follows from \cite[Theorem 1]{harris_optimal_2021} that the convex relaxation \Cref{prob:CP_Mayer} is lossless with respect to \Cref{prob:MICP_Mayer}. Consequently, $\big(u_\star(t),\mu_\star(t)\big)$ is also an optimal control of \Cref{prob:MICP_Mayer} and $\big(u_\star(t),\mu_\star(t)\big)\in\mathcal{U}_e$ for almost every $t\in[0,t_f]$. Moreover, it follows from \Cref{lemma:vertices_normality} that $\big(u_\star(t),\mu_\star(t)\big)$ is piecewise constant on $[0,t_f]$.
Since $\big(u_\star(t),\mu_\star(t)\big)\in\mathcal{U}_e$, it follows from \cref{eq:def_augmented_set} that $u_\star(t)\in\SU$ for almost every $t\in[0,t_f]$. 

Since $\mathcal{U}_e\subseteq\tilde{\mathcal{U}}_e$ and $\hat{u}^*(t)$ minimizes $H_2$ over $\tilde{\mathcal{U}}_e$, it also minimizes $H_2$ over $\mathcal{U}_e$. Hence, $\hat{u}^*(t)$ satisfies \cref{eq:optimal_control_p2}. Applying \Cref{lemma:MICP2mayer}, we conclude that $\mu_\star(t)=\psi(u_\star(t))$ and that $u_\star(t)$ satisfies \cref{eq:optimal_control_p1}.
\end{proof}

\Cref{thm:CP_for_MICP} formalizes the equivalence between the relaxed convex program \Cref{prob:CP_Mayer} and the original mixed-integer problem \Cref{prob:MICP}, establishing that this formulation constitutes a \textit{lossless convexification}. By mapping the problem to a standard convex program, the optimization can be solved in polynomial time using interior-point methods or first-order solvers \cite{boyd_convex_2023}. Crucially, because the convexification is lossless, the optimal trajectory of the relaxed convex program satisfies the original discrete control constraints without requiring any ad-hoc heuristic rounding or post-processing. This guarantees both exactness and algorithmic reliability, unlocking the potential for executing complex mixed-integer optimal control strategies on resource-constrained embedded systems, such as onboard spacecraft guidance.

\section{Numerical Simulations}\label{sec:NumericalSim}

We illustrate the application of this lossless convexification framework through a fuel-optimal problem ($\psi(u) = \|u\|_1$) of a low-Earth orbit (LEO) satellite rendezvous maneuver in an elliptical orbit. Consider a non-cooperative chief spacecraft and a chaser spacecraft with cold-gas thruster actuators in its three orthogonal directions, which are saturated at maximum acceleration $u_{\max} = 0.05 \frac{\rm m}{\rm s^2}$, a constraint typical in small satellites \cite{weston_state---art_2025}. The input set is the same as in \Cref{ex:set_example_1} and is given by
\begin{equation*}
    \mathcal{U} \triangleq \{-u_{\max},0,u_{\max}\}^3.
\end{equation*}
The relative motion of a chaser spacecraft with respect to a chief in an elliptic orbit can be described by Yamanaka--Ankersen (YA) equations \cite{yamanaka_new_2002}. The relative state is $x(t) \triangleq \begin{bmatrix} r(t)^\top & v(t)^\top \end{bmatrix}^\top$, where $r(t) = \begin{bmatrix} x(t) & y(t) & z(t) \end{bmatrix}^\top$ are the radial, along-track, and cross-track displacements, respectively, and $v \triangleq \dot{r}$. 
The system matrices in \cref{eq:linear_system} are given by
\begin{equation}
    A(t) =
    \begin{bmatrix}
    0 & 0 & 0 & 1 & 0 & 0 \\
    0 & 0 & 0 & 0 & 1 & 0 \\
    0 & 0 & 0 & 0 & 0 & 1 \\
    \dot{\theta}(t)^2 + \dfrac{2\mu}{r_c(t)^3}
        & \ddot{\theta}(t)
        & 0
        & 0
        & 2\dot{\theta}(t)
        & 0 \\
    -\ddot{\theta}(t)
        & \dot{\theta}(t)^2 - \dfrac{\mu}{r_c(t)^3}
        & 0
        & -2\dot{\theta}(t)
        & 0
        & 0 \\
    0
        & 0
        & -\dfrac{\mu}{r_c(t)^3}
        & 0
        & 0
        & 0
    \end{bmatrix},
    \qquad
    B=
    \begin{bmatrix}
    0_{3\times 3} \\
    I_{3\times 3}
    \end{bmatrix},
    \label{eq:ya_system_matrices}
\end{equation}
where $r_c(t) \in \BBR$ is the instantaneous orbital radius of the chief,
$\theta(t) \in \BBR$ is its true anomaly, and $\mu$ is the gravitational parameter
of the central body. For a Keplerian chief orbit with semimajor axis
$a>0$, eccentricity $e \in [0,1)$ and specific angular momentum $h \triangleq \sqrt{\mu a(1-e^2)}$, the time-varying orbital parameters in \cref{eq:ya_system_matrices} are given by
\begin{equation*}
    \dot{r}_c(t)
    =
    \frac{he}{a(1-e^2)}\sin\theta(t), \quad
    r_c(t)
    =
    \frac{a(1-e^2)}{1+e\cos\theta(t)},
    \quad
    \ddot{\theta}(t)
    =
    -\frac{2\mu e\sin\theta(t)
    \left(1+e\cos\theta(t)\right)^3}{(a(1-e^2))^3}, \quad 
    \dot{\theta}(t)
    =
    \frac{h}{r_c(t)^2}.
\end{equation*}
Throughout this section we consider a chief spacecraft in an elliptical orbit with an eccentricity of $e=0.7$ and a semimajor axis of $a=7102.8 \times 10^3 \text{ m}$.

\subsection{Algorithm Implementation and Time-Horizon Selection}

% To solve \Cref{prob:CP_Mayer} numerically, we discretize the time horizon $[0,t_f]$ into $N$ equal intervals of duration $\Delta t = t_f / N$, defining the discrete time instances $t_k = k \Delta t$ for $k = 0, 1, \dots, N$. Applying a zero-order hold (ZOH) to the control input over each interval, the discrete-time matrices are given by
% %
% \begin{equation*}
%     \hat{A}_k = \exp\left(\hat{A}(t_k)\Delta t\right), \quad \hat{B}_k = \left(\int_{0}^{\Delta t} \exp\left(\hat{A}(t_k)\tau\right)\mathrm{d}\tau\right)\hat{B}(t_k),
% \end{equation*}
% %
% where $\hat{x}_k = [x_k^T, x_{c,k}]^T$ and $\hat{u}_k$ denote the augmented state and control vectors at the $k$-th sampling node, respectively. This reformulation results in the following optimization problem:
To solve \Cref{prob:CP_Mayer} numerically, we discretize the time horizon $[0,t_f]$ into $N$ equal intervals of duration $\Delta t = t_f / N$, defining the discrete time instances $t_k = k \Delta t$ for $k = 0, 1, \dots, N$. Applying a zero-order hold (ZOH) to the control input over each interval, the discrete-time matrices are given by $\hat{A}_k = \exp\left(\hat{A}(t_k)\Delta t\right), \quad \hat{B}_k = \left(\int_{0}^{\Delta t} \exp\left(\hat{A}(t_k)\tau\right)\mathrm{d}\tau\right)\hat{B}(t_k)$, where $\hat{x}_k = [x_k^T, x_{c,k}]^T$ and $\hat{u}_k$ denote the augmented state and control vectors at the $k$-th sampling node, respectively. This reformulation results in the following optimization problem:
\begin{problem}\rm{Discrete-time Mayer-form Convex Problem}\label{prob:CP_Mayer_discrete}
    \begin{subequations}
        \begin{align}
            \min \quad & x_{c,N} \label{eq:Discrete_Objective}\\
            \text{subject to} \quad 
            & \hat{x}_{k+1} = \hat{A}_k \hat{x}_k + \hat{B}_k \hat{u}_k, \quad k = 0, \dots, N-1, \label{eq:Prob_Discrete_Dynamics}\\
            & \hat{u}_k \in \tilde{\mathcal{U}}_e, \quad k = 0, \dots, N-1, \label{eq:Prob_Discrete_Control}\\
            & x_0 = x_0, \quad x_N = 0, \label{eq:Prob_Discrete_BC}\\
            & x_{c,0} = 0. \label{eq:Prob_Discrete_Cost_Init}
        \end{align}
    \end{subequations}
\end{problem}

This discrete finite-dimensional problem is solved using CVXPY~\cite{diamond_cvxpy_2016} with the ECOS solver on a MacBook Air equipped with an Apple M1 chip and 16~GB of RAM.\footnote{The implementation is available at: \url{https://github.com/FelipeArenasUribe/MICP_LosslessConvexification}}

In order to solve \Cref{prob:CP_Mayer_discrete}, one must provide a feasible time horizon $t_f$. Selecting an inappropriate $t_f$ can result in an infeasible optimization problem if the horizon is too short to satisfy the boundary conditions under the control bounds. Alternatively, a long time-horizon results in an undesirable maneuver for mission constraints. To solve this challenge, \Cref{prob:CP_Mayer_discrete} can be embedded within a one-dimensional outer-loop optimization algorithm. Let $0<\underline{t}_f < \bar{t}_f<\infty$ be user-defined time limits and consider the time-horizon interval $[\underline{t}_f, \bar{t}_f]$. To find a feasible time-horizon $t^*_f \in [\underline{t}_f, \bar{t}_f]$, a bisection method can be utilized over the interval to locate the smallest $t_f$ for which \Cref{prob:CP_Mayer_discrete} returns a feasible solution. 

As illustrated in \Cref{fig:bisection}, the algorithm initializes the dynamic iteration bounds as $t_\ell = \underline{t}_f$ and $t_u = \bar{t}_f$. At each iteration, the optimization problem is solved using the midpoint horizon $t_{\mathrm{mid}} = (t_\ell + t_u)/2$. If the problem is feasible, the upper bound is updated ($t_u \leftarrow t_{\mathrm{mid}}$); otherwise, the lower bound is updated ($t_\ell \leftarrow t_{\mathrm{mid}}$). This process repeats until the interval narrows below a user-defined convergence tolerance ($t_u - t_\ell \leq \varepsilon$), yielding the optimal feasible horizon $t_f^*$. The bisection algorithm requires a feasible upper bound $\bar{t}_f$, while the lower bound $\underline{t}_f$ may be infeasible.

\begin{figure}[ht]
\centering
\begin{tikzpicture}[
    node distance = 5mm and 6mm, % Reduced from 6mm and 10mm
    box/.style     = {rectangle, rounded corners=3pt, draw, minimum width=22mm, % Narrowed width
                      minimum height=8mm, align=center, font=\small}, % Lowered height
    decision/.style= {diamond, draw, aspect=1.5, align=center, font=\small,
                      inner sep=2pt},
    arr/.style     = {-{Stealth[length=4pt]}, thin},
    labl/.style    = {font=\scriptsize, inner sep=1pt},
]

%% Nodes (left to right)
\node[box] (init)
    {Initialize\\$t_\ell = \underline{t}_f,\; t_u = \bar{t}_f$};

\node[decision, right=of init] (conv)
    {$t_u - t_\ell \leq \varepsilon$?}; % Combined onto one line to save height

\node[box, right=of conv] (mid)
    {$t_{\mathrm{mid}} = \tfrac{t_\ell+t_u}{2}$};

\node[box, right=of mid] (solve)
    {Solve \Cref{prob:CP_Mayer_discrete}\\$t_f = t_{\mathrm{mid}}$};

\node[decision, right=of solve] (feas)
    {Feasible?};

%% Convergence: YES now exits downward/sideways to keep the profile low
\node[box, below=8mm of conv] (done)
    {$t_f^{\ast} = t_u$};

%% Main flow arrows
\draw[arr] (init)  -- (conv);
\draw[arr] (conv)  -- node[labl,above]{No}  (mid);
\draw[arr] (mid)   -- (solve);
\draw[arr] (solve) -- (feas);

%% Converged: exit downward (saves massive vertical headspace)
\draw[arr] (conv) -- node[labl,right]{Yes} (done);

%% Feasible → tighten upper bound → loop back (Reduced overshoot to 6mm)
\draw[arr] (feas.north) -- ++(0, 6mm)
    node[labl,right,pos=0.5]{Yes, $t_u \leftarrow t_{\mathrm{mid}}$}
    -| (conv.north);

%% Infeasible → raise lower bound → loop back (Reduced overshoot to 6mm)
\draw[arr] (feas.south) -- ++(0,-6mm)
    node[labl,right,pos=0.5]{No, $t_\ell \leftarrow t_{\mathrm{mid}}$}
    -| (conv.south);

\end{tikzpicture}
\caption{Outer bisection loop for time horizon selection.}
\label{fig:bisection}
\end{figure}

The output of \Cref{prob:CP_Mayer_discrete} is a control sequence $\{u_k\}_{k=0}^{N-1}$. To evaluate the quality of the solutions, we introduce the following metrics: distance of the input to the discrete set and the average distance of the control sequence to the discrete set, which are defined as
\begin{equation*}
    d(u_k) \triangleq \min_{\zeta \in \SU} \| u_k - \zeta \|_2, \quad 
    \bar{d}(u) \triangleq \frac{1}{N-1} \sum_{k=0}^{N-1} \frac{d(u_k)}{u_{\max}}.
\end{equation*}
These metrics characterize the discreteness of the solution. Lower values in both metrics imply the control trajectory is closer to satisfying the input constraint $u \in \SU$.

\subsection{Validation of Lossless Convexification}\label{subsec:validation}

We first demonstrate that the resulting optimal control sequence from \Cref{prob:CP_Mayer_discrete} is discrete-valued and compare it with an approximate solution of \Cref{prob:MICP}.
The initial conditions are $r_0 = [\begin{smallmatrix} -100 & -500 & -100\end{smallmatrix}]\,\text{m}$ and $v_0 = [\begin{smallmatrix} 0 & 0 & 0 \end{smallmatrix}]\,\text{m/s}$. 
The optimal time horizon is $t_f^* = 192.1875\,\text{ s}$, which was computed using the bisection algorithm from \Cref{fig:bisection} with $\bar{t}_f = 500.0\,\text{ s}$, $\underline{t}_f=100.0\,\text{ s}$, $\varepsilon=1.0\,\text{ s}$ and a discretization grid of $N=200$. The bisection loop converges to a solution in 0.4458 seconds. 

To closely approximate the continuous-time solution of \Cref{prob:CP_Mayer}, we solve \Cref{prob:CP_Mayer_discrete} with a discretization grid of $N=1000$. For real-time applications, coarser grids may be adopted to reduce computational effort, as discussed in \cref{subsec:preserving_disc}.
We compute the solution of \Cref{prob:MICP_Mayer} as a finite-dimensional mixed-integer linear program (MILP) using Gurobi. Due to software constraints, the mixed-integer problem uses $N=200$.
Additionally, we relaxed the terminal constraint with tolerances of $0.5\,\text{m}$ in position and $0.005\,\text{m/s}$ in velocity to facilitate convergence to a solution.

\begin{figure}[ht]
    \centering
    % First Subfigure (Left)
    \begin{subfigure}[b]{0.48\linewidth}
        \centering
        \includegraphics[width=0.8\linewidth]{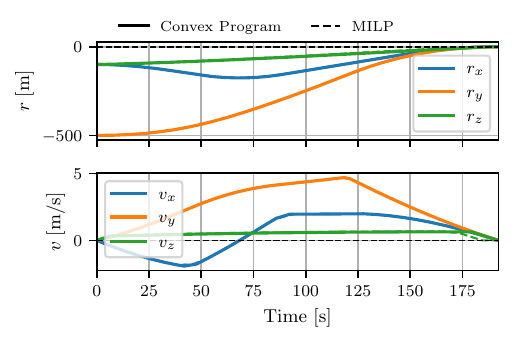}
        \caption{System states during rendezvous, both solution trajectories show equal convergence to the desired terminal state.}
        \label{fig:RendezvousStates}
    \end{subfigure}
    \hfill % Pushes the subfigures to the left and right edges
    % Second Subfigure (Right)
    \begin{subfigure}[b]{0.48\linewidth}
        \centering
        \includegraphics[width=0.9\linewidth]{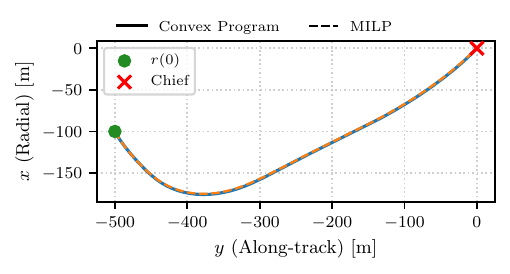}
        \caption{Trajectory of fuel-optimal rendezvous maneuver in the radial and along-track plane.}
        \label{fig:Rendezvoustrajectory}
    \end{subfigure}
    
    \caption{Spacecraft states and relative 2D trajectory during the fuel-optimal rendezvous maneuver for the solutions obtained with the CP and the MILP.}
    \label{fig:RendezvousCombined} % Main label for the entire group
\end{figure}

The MILP optimal cost is $x_{c,N} = 19.1227$, while the convex-program optimal cost is $x_{c,N} = 19.3827$. 
The MILP objective is lower due to the relaxed terminal constraint.
The computation time of \Cref{prob:CP_Mayer_discrete} is $0.9782$~seconds, while the computation time of the MILP is $473.18$~seconds. As expected, the CP was approximately $484$ times faster.
\Cref{fig:RendezvousStates} shows the solution of the system's states. Moreover, \Cref{fig:Rendezvoustrajectory} displays the in-plane trajectory of the spacecraft in the Hill frame.

\begin{figure}[ht]
    \centering
    \includegraphics[width=0.38\linewidth]{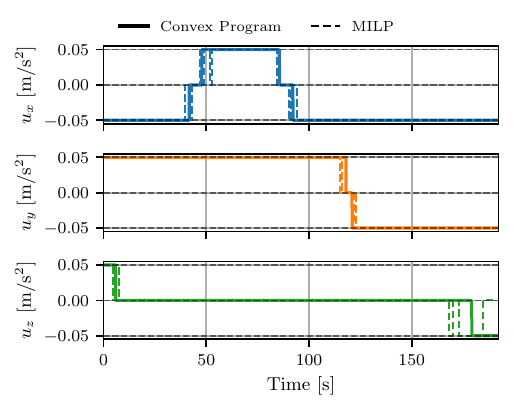}
    %\vspace{-0.25in}
    \caption{Optimal control inputs with respect to the admissible discrete set $\SU$.}
    \label{fig:RendezvousInputs}
    %\vspace{-0.1 in}
\end{figure}

The control inputs are shown in \Cref{fig:RendezvousInputs}. 
It can be seen that for all $t \in [0,t_f]$, the optimal control satisfies the input constraint $u(t) \in \SU$. 
The average distance to the discrete set is $\bar{d}(u) = 0.002148$. 
This indicates that the control sequence remains effectively discrete-valued, despite being computed via a finite-dimensional approximation of \Cref{prob:MICP}. Finally, note that the control sequence for the CP has a similar profile to that of the MILP.

\subsection{Preservation of Discrete-Valued Controls Under Reduced Grid Resolution}\label{subsec:preserving_disc}

Real-time implementation requires balancing solver execution time with the need for a discrete-valued control. This introduces an inherent trade-off when selecting the number of sampling instances. While increasing $N$ aligns the finite-dimensional approximation \Cref{prob:CP_Mayer_discrete} closely with \Cref{prob:CP_Mayer}, it simultaneously drives up solver time. Conversely, reducing $N$ expedites the solver at the cost of degrading the discrete-valued nature of the control trajectory, thereby weakening its adherence to the input constraint \cref{eq:input_cons} and the optimal solution of \Cref{prob:CP_Mayer}.

To evaluate this trade-off we conduct a parametric study analyzing solution discreteness and solver efficiency across a range of grid sizes. For a fixed time horizon $t_f = 200\,\text{s}$, a batch of 10 initial conditions is uniformly sampled from the set $\SX_0 \triangleq \{x \in \BBR^6 : \|r\|_{\infty} \leq 500, \|v\|_{\infty} \leq 5\}$. For each initial condition, \Cref{prob:CP_Mayer_discrete} is solved using a different quantity of sampling grid points, sweeping from an initial value of $N_0 = 100$ to a final value of $N_f = 1000$ in increments of 100.

\begin{figure}[ht]
    \centering
    \includegraphics[width=0.43\linewidth]{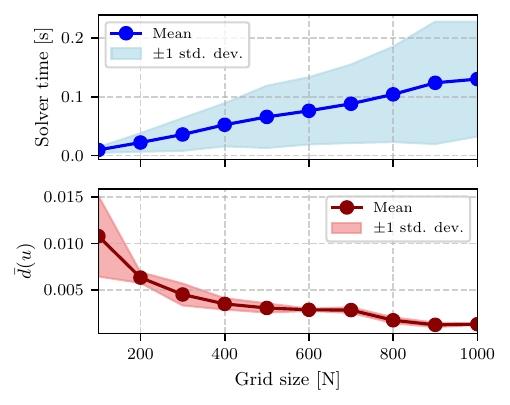}
    \vspace{-0.1in}
    \caption{Effect of grid size on solver time and control discreteness. Increasing the number of grid points improves adherence to the discrete control set but increases computational time.}
    \label{fig:Grid_Tradeoff}
\end{figure}

The numerical results of the parametric sweep are illustrated in \Cref{fig:Grid_Tradeoff}, which details the mean solver time and the mean normalized average distance $\bar{d}(u)$ from the different initial conditions across the evaluated grid sizes. 
The mean solver time exhibits a strictly linear dependence on the number of grid points $N$, scaling from $0.025\,\text{seconds}$ at $N=200$ to $0.135\,\text{seconds}$ at $N=1000$. 
It can be seen that solver time varies based on the initial condition, and this variance increases as the grid size increases. 
Conversely, the mean normalized average distance decays approximately exponentially as grid density increases. 
Increasing $N$ beyond a certain threshold yields negligible improvements in control discreteness at the expense of increased computation time. Ultimately, the remaining quantization error at high grid sizes is small enough to be readily compensated for through low-level feedback control, validating the framework's suitability for real-time applications. Moreover, unlike solver time, the variance reduces as the grid size increases, which demonstrates that control discreteness does not depend on initial conditions but is strictly a function of the grid size. 

\subsection{Real-Time Performance}\label{subsec:performance}

Real-time spacecraft control guidelines from NASA's Safe \& Precise Landing Integrated Capabilities Evolution (SPLICE) project specify an attitude and position guidance update-rate requirement of 3 seconds, with a target goal of 1 second \cite{sostaric_splice_2021}. 
Although we focus only on position control, this benchmark provides a meaningful reference for evaluating the computational performance of our algorithm. To assess real-time feasibility, we perform a Monte Carlo simulation in which 1000 initial conditions are uniformly sampled from $\SX_0$. 
For each sampled initial condition, \Cref{prob:CP_Mayer_discrete} is solved with a time horizon of $t_f = 400\,\text{s}$ and $N = 400$ grid points based on the results from \cref{subsec:preserving_disc}.

\begin{figure}[ht]
    \centering
    % First Subfigure (Left)
    \begin{subfigure}[b]{0.48\linewidth}
        \centering
        \includegraphics[width=0.9\linewidth]{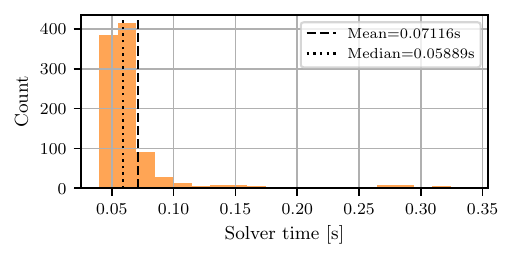}
        \caption{Histogram of solver times. All solver times are below 1\,s, demonstrating real-time feasibility for the selected grid.}
        \label{fig:Solvertime_hist}
    \end{subfigure}
    \hfill % Distributes the subfigures evenly across the text width
    % Second Subfigure (Right)
    \begin{subfigure}[b]{0.48\linewidth}
        \centering
        \includegraphics[width=0.9\linewidth]{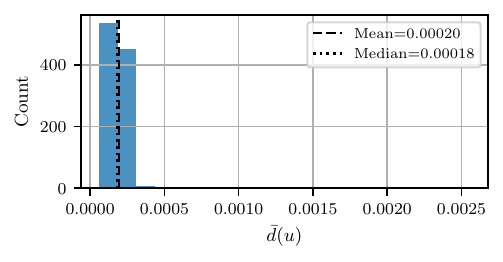}
        \caption{Histogram of the average distance to the discrete set $\bar{d}(u)$. The algorithm reliably produces discrete-valued controls.}
        \label{fig:Distance_hist}
    \end{subfigure}
    
    \caption{Monte Carlo simulation results showing distribution of solver execution times and the distance metrics to the discrete control set.}
    \label{fig:MonteCarloHistograms}
\end{figure}

\begin{figure}[ht]
    \centering
    \includegraphics[width=0.5\linewidth]{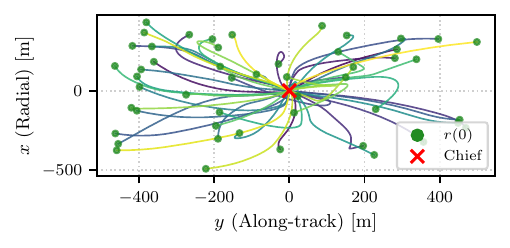}
    \caption{Subset of 50 trajectories from the 1000 trajectories computed in the Monte Carlo simulation.}
    \label{fig:mc_trajs}
\end{figure}

The Monte Carlo simulation demonstrates that the solver time across all evaluated initial conditions remains well below the 1-second threshold (see \Cref{fig:Solvertime_hist}), satisfying the SPLICE target, a subset of the computed trajectories is displayed in \Cref{fig:mc_trajs}. 
The resulting profile exhibits a right-skewed distribution, featuring a mean solver time of $0.07116$~seconds and a median of $0.05889$~seconds. This indicates that while a small subset of initial conditions creates a tail extending toward $0.40$~seconds, the majority of cases resolve substantially faster than the target requirement. Furthermore, as shown in \Cref{fig:Distance_hist}, the distribution of the normalized average distance to the discrete control set, $\bar{d}(u)$, is narrowly clustered around a mean of approximately $0.0002$. This highly concentrated profile demonstrates that the algorithm consistently computes near-discrete control inputs across the entire operational envelope, which aligns with the results in \cref{subsec:preserving_disc}. Overall, these results confirm that the proposed method simultaneously guarantees real-time performance and control discreteness.

\section{Conclusions}\label{sec:Conclusions}

This paper presents a lossless convexification for a class of optimal control problems of linear time-varying systems with discrete-valued control inputs. 
This result is anchored by two primary theoretical contributions. First, \Cref{lemma:Normality_augmented} establishes that system normality is preserved when transforming between the Lagrange and Mayer formulations using an epigraph transformation.
Next, by extending the foundations from \cite{harris_optimal_2021}, \Cref{thm:CP_for_MICP} proves that using a convex hull relaxation of the input constraint provides a lossless convexification of \Cref{prob:MICP}.
This formulation allows the mixed-integer problem to be solved reliably as a convex program. To bridge the gap between theoretical guarantees and algorithmic deployment, we introduced a finite-dimensional approximation alongside an iterative bisection algorithm designed to identify a feasible time horizon and compute the optimal control sequence.

Numerical simulations demonstrate that the finite-dimensional algorithm consistently yields near-discrete controls while remaining within real-time computational constraints. 
These results confirm the practical feasibility of the proposed framework for safety-critical aerospace applications. 
Future work will extend the framework to incorporate state constraints while preserving input discreteness. 
Another direction is to integrate the algorithm as a guidance layer with a feedback control layer that compensates for approximation errors while enforcing input and state constraints. 
This integration could extend the open-loop lossless convexification framework into a closed-loop architecture for real-time control with discrete-valued inputs.

\section*{Acknowledgments}

The authors acknowledge the use of large language models to assist with the development of the numerical implementation. All theoretical developments and scientific conclusions are original contributions of the authors.

\bibliography{References}

\end{document}